\documentclass[12pt]{article}
\usepackage{mystyle}

\usepackage{parskip}
\usepackage{mathtools}
\usepackage{stmaryrd}

\usepackage[
backend=biber,
style=numeric,
sorting=nyt
]{biblatex}
\usepackage{cleveref}
\newcommand{\simp}[2]{V_{\overline{#1}, #2}}

\title{Fusion Rings over Drinfeld Doubles}
\author{Wenqi Li \\ wl2935@columbia.edu}
\date{\today}

\makeatletter
\let\runauthor\@author
\makeatother

\begin{document}

\maketitle

\begin{abstract}

The fusion rules in $\mathrm{Rep}_f D(G)$ for a finite group $G$ can be computed in terms of character inner products. Using an explicit formula (Theorem \ref{actionfusion}) for these fusion rules, we show that $\mathrm{Rep}_f D(G)$ is multiplicity free for two infinite families of finite groups: the Dihedral groups and the Dicyclic groups. In fact, we will compute all fusion rules in these categories. Multiplicity freeness is a desired property for modular tensor categories, since it greatly simplifies the computation of $F$-matrices. Furthermore, we observe that the fusion rules for Dihedral groups $D_{2n}$ with $n$ odd are extremely similar to the fusion rules of Type $B$ level $2$ fusion algebras of Wess-Zumino-Witten conformal field theories. Moreover, we give a proof of the fusion rule formula by using Mackey theory.
\end{abstract}

\section{Introduction}

A modular tensor category is a semisimple ribbon category with finitely many simple objects, with an invertible matrix called the $S$-matrix. %whose rows and columns correspond to simple objects. entries are certain scalar maps associated with the simple objects. 
The properties of the $S$-matrix can be used to give a projective representation of $\SL_2(\Z)$, hence the name ``modular''. Modular tensor categories have important applications in knot theory and topological quantum computing, since in some sense (\cite{Turaev+2010}) each modular tensor category uniquely determines a topological quantum field theory. 

One important source of examples of modular tensor categories comes from a construction called the \textit{Drinfeld double}. The Drinfeld double $D(H)$ of a Hopf algebra $H$ is a key construction in various aspects of representation theory and mathematical physics. An important class of concrete examples is the Drinfeld double of a finite group $G$, where the Hopf algebra $H$ is the group algebra $kG$ for $k$ a field. It is customary to denote such a Drinfeld double by $D(G)$ rather than $D(kG)$. The categories of finite dimension representations of Drinfeld doubles of finite groups, $\mathrm{Rep}_f D(G)$, are modular tensor categories. As a generalization, \cite{DIJKGRAAF199160} constructed and developed the theory of twisted Drinfeld doubles of finite groups, where the Drinfeld double is twisted by a 3-cocycle $\omega \in Z^3(G, \C^{\times})$. In \cite{mason_campbell_robertson_hurley_tobin_ward_1995}, Mason explains the application of Drinfeld doubles of finite groups to conformal field theory, specifically to holomorphic vertex operator algebras and holomorphic orbifolds. The twisted Drinfeld double was generalized even further by \cite{GOFF20103007} by replacing the group algebra $kG$ by the group algebra of some quotient of $G$. The twisted Drinfeld doubles of finite groups have been of great interest in research on modular tensor categories and holomorphic vertex operator algebras. (See for example \cite{Mason_Ng_twist} for conjectures in the case of rational orbifolds and many examples.)

Computing the fusion rules, $S$-matrices, $R$-matrices and $F$-matrices are the central tasks in studying modular tensor categories. In the case of $\mathrm{Rep}_f D(G)$, the $S$-matrix is always known, and it is of interest to compute the $R$-matrix and the $F$-matrices for various examples. The $F$-matrices are of interest, since, for example, they are needed to make braid group representations explicit (see \cite{Trebst2008ASI} for a discussion). \cite{RSW_classification} provides a classification of modular tensor categories of small ranks, one of which is realized as $\mathrm{Rep}_f D(G)$ for some group $G$, and their $F$-matrices are computed. The modular tensor categories listed in \cite{RSW_classification} of the form $\mathrm{Rep}_f D(G)$ has trivial $F$-matrices, but this is not true in general. For instance, the paper \cite{Cui_Hong_Wang_DS3} computes the non-trivial $F$-matrix for $D(S_3)$.

For fusion rules, \cite{GoffAbExt} obtains explicit formulas in the general case of the representation category of a twisted Drinfeld double of a finite group $G$, which specializes to a formula for the fusion rules in the untwisted representation category $\mathrm{Rep}_f D(G)$. In this paper, we give a proof of the formula for fusion rules in $\mathrm{Rep}_f D(G)$ using Mackey theory. 

We say a fusion category is \textit{multiplicity free} if all the fusion coefficients $N_{ij}^k$ are bounded above by $1$. Multiplicity freeness is a desirable property, since it greatly simplifies the calculation of $F$-matrices. As far as we know, there are no known cases where the $F$-matrices are computed in a category that is not multiplicity free. Therefore, it is of interest to know which finite groups give Drinfeld doubles whose module categories are multiplicity free.

The main result of this paper is the following: we show that $\mathrm{Rep}_f D(G)$ is multiplicity free for two infinite families of finite groups: the Dihedral groups and the Dicyclic groups. In fact, we will compute all fusion rules in these categories. Furthermore, we observe that the fusion rules for Dihedral groups $D_{2n}$ with $n$ odd is extremely similar to the fusion rules of Type $B$ level $2$ fusion algebras of Wess-Zumino-Witten conformal field theories.

\section{Modular Tensor Categories from Drinfeld Doubles}\label{DR}

Let $G$ be a finite group. We list here the structure of the \textit{Drinfeld double} $D(G)$ as a Hopf algebra. As a vector space it is $\C[G] \otimes \cur{O}(G)$ where $\cur{O}(G)$ denotes the algebra of functions on $G$, and its Hopf algebra structure is given by the following:
\begin{align*}
\text{multiplication} & \quad (x \otimes \delta_g)(y \otimes \delta_h) = \delta_{gx, xh}(xy \otimes \delta_g) \\
\text{unit} & \quad 1 = \sum_{g \in G} e \otimes \delta_g \\
\text{comultiplication} & \quad \Delta(x \otimes \delta_g) = \sum_{g_1 g_2 = g} (x \otimes \delta_{g_1}) \otimes (x \otimes \delta_{g2}) \\
\text{counit} & \quad \epsilon(x \otimes \delta_g) = \delta_{g, e} \\
\text{antipode} & \quad \gamma(x \otimes \delta_g) = x^{-1} \otimes \delta_{x^{-1}g^{-1}x}
\end{align*}

We investigate the category of finite dimensional representations of $D(G)$. Denote this category by $\mathrm{Rep}_f(D(G))$. 

Here is a description of the simple objects as in \cite{gould_1993}. Let $K \subset G$ be a conjugacy class, $g_K$ be a fixed representative of $K$, and $\pi$ an irreducible representation of $C(g_K)$, the centralizer of $g_K$. For each $s \in K$, fix a $\tau_s$ such that $\tau_s g_K \tau_s^{-1} = s$. This is the same as choosing a complete set of coset representatives in $G/C(g_K)$. 

Let $V_{\pi}$ be the $\C[C(g_K)]$-module corresponding to the representation $\pi$. Define
\[
V_{g_K, \pi} = \C[G] \otimes_{\C[C(g_K)]} V_{\pi}.
\]
So as a $\C[G]$ module, $V_{g_K, \pi}$ is the induced representation of $V$ to $G$. Since $\tau_s$'s are a complete set of coset representative, we know that
\[V_{g_K, \pi} \cong \bigoplus_{s \in K} \tau_s \otimes V_{\pi},\]
with a basis of the form $\{\tau_s \otimes v_i \}_{s \in K}$ if $v_i$'s are a basis of $V_{\pi}$. 

The $G$ action is given by $g(\tau_s \otimes v) = \tau_{gsg^{-1}} \otimes (\tau_{gsg^{-1}}^{-1} g \tau_s)v$ on the basis elements. This action can be extended to a $D(G)$-action by defining
\[
\delta_f (\tau_s \otimes v) = \delta_{f, s}(\tau_s \otimes v).
\]

It is well-known that these representations are indeed simple, and they are all the simple objects (see \cite{gould_1993} for a proof).

The category $\mathrm{Rep}_f(D(G))$ is modular tensor, and its $S$-matrix can be explicitly computed.

\begin{Theorem}\label{DGMTC}
The category of finite dimensional representations of $D(G)$ is a modular tensor category. For $g \in G$, denote its conjugacy class by $\overline{g}$. The dual of a simple object is
\begin{equation}\label{DGdual}
\simp{g}{\pi}^* \cong \simp{g^{-1}}{\pi^*}
\end{equation}
and the normalized $S$-matrix of this category is
\begin{equation}\label{DGs}
s_{(\overline{g}, \pi), (\overline{g'}, \pi)} = \frac{1}{|C(g)|}\frac{1}{|C(g')|} \sum_{\substack{h \in G, \\ hg'h^{-1} \in C(g)}} \tr_{\pi}(hg'h^{-1}) \tr_{\pi'}(h^{-1}g^{-1}h).
\end{equation}
\end{Theorem}
\begin{proof}
This is Theorem 3.2.1 in \cite{bakalov2001lectures}.
\end{proof}

\section{A Character Formula for Fusion Rules}

In this section, we introduce a formula (Theorem \ref{generalfusion}) for the fusion rules in the category $\mathrm{Rep}_f D(G)$ in terms of character inner products. It is true that the Verlinde formula gives us the fusion rules, but it is hard to compute the summation in the Verlinde formula because it involves the $S$-matrix. The formula in Theorem \ref{generalfusion} only requires certain group theoretic data of $G$.

We have seen above that irreducible representations $V_{\overline{g}, \pi}$ of the Drinfeld double $D(G)$ are determined by a conjugacy $\overline{g}$ for some $g \in G$ and an (isomorphism class of) irreducible representation $\pi$ of the centralizer $C(g)$. Therefore, we can use the character $\chi$ of the representation $\pi$ to compute the character $\widehat{\chi}_{\overline{g}}$ of $V_{\overline{g}, \pi}$.

\begin{Proposition}\label{charformula}
Let $G$ be a finite group and let $K$ be a conjugacy class in $G$. Let $V_{K, \pi} = V_{\overline{g_K}, \pi}$ be an irreducible representation of $D(G)$, where $g_K$ is a chosen representative of $K$ and $\pi$ is an irreducible representation of $C(g_K)$. If $\chi$ is the character of $\pi$, then the character $\widehat{\chi}_K$ of $V_{K, \pi}$ is given by
\[
\widehat{\chi}_K(x \otimes \delta_g) = \delta_{g \in K}\delta_{x \in C(g)} \chi(rxr^{-1}),
\]
where $r$ is such that $rgr^{-1} = g_K$.
\end{Proposition}
\begin{proof}
This is Proposition 2.2.4 in \cite{fusionrule}.
\end{proof}

We will use the above proposition to prove the following: 

\begin{Theorem}\label{actionfusion}
Let $G$ be any finite group. Let $K, L, J$ be conjugacy classes of $G$, and fix representatives $g_K$, $g_L$ and $g_J$. Let $C(g)$ denote the centralizer of $g \in G$ as usual. Let $\chi, \psi, \varphi$ be characters of $C(g_K), C(g_L), C(g_J)$. Then
\begin{equation}\label{actionformula}
    \langle \widehat{\chi} \otimes \widehat{\psi}, \widehat{\varphi} \rangle = \sum_{(k, l) \in O} \langle\chi^{(k)} \psi^{(l)}, \varphi\rangle_{C(k) \cap C(l)}.
\end{equation}
where $O = \{(k, l) \in K \times L \mid kl = g_J\}/\sim$ and $(k, l)\sim (k', l')$ if there exists some $\gamma \in C(g_J)$ such that $(k', l') = (\gamma k \gamma^{-1}, \gamma l \gamma^{-1})$.
\end{Theorem}

Some explanations are in order. First note that $\chi^{(k)}$ (same for $\psi^{(l)}$) makes sense and is well-defined. First of all, if $r^{-1}g_Kr = k$, then for $x \in C(k)$, $rxr^{-1}$ centralizes $g_K$ since
\[rxr^{-1}g_Krx^{-1}r^{-1} = rxkx^{-1}r^{-1} = rkr^{-1} = g_K.\]
Hence $\chi^{(k)}$ is defined on $C(k)$. Moreover, the definition of $\chi^{(k)}$ does not depend on the choice of $r$. If $r_1$ and $r_2$ both satisfy $r_1^{-1}g_K r_1 = k = r_2^{-1} g_K r_2$, then $r_2r_1^{-1}g_K r_1 r_2^{-1} = g_K$, so $r_1r_2^{-1}$ is in $C(g_K)$. Since $\chi$ is a class function on $C(g_K)$, we see that $\chi(x) = \chi(r_2r_1^{-1} x r_1r_2^{-1})$ for all $x \in C(g_K)$. Letting $y = r_1^{-1}xr_1$, we obtain that $\chi(r_1yr_1^{-1}) = \chi(r_2yr_2^{-1})$ for all $y \in C(k)$ since conjugation by $r_1$ is clearly a bijection between $C(g_K)$ and $C(k)$.

\cite{fusionrule} proves this theorem in the special case where all centralizers are assumed to be normal. Theorem \ref{actionfusion} here is a generalization. In \cite{GOFF20103007} and \cite{GoffAbExt}, a more general formula is proved for twisted Drinfeld doubles, which would imply Theorem \ref{actionfusion} in the case where the twisting is trivial. 

An equivalent formulation is the following:

\begin{Theorem}\label{generalfusion}

\begin{equation}\label{prettyformula}
\langle \widehat{\chi} \otimes \widehat{\psi}, \widehat{\varphi} \rangle = \frac{|J|}{|G|} \sum_{\substack{(k, l) \in K \times L,\\kl = g_J}} |C(k) \cap C(l)| \langle\chi^{(k)} \psi^{(l)}, \varphi\rangle_{C(k) \cap C(l)}.
\end{equation}
where $\chi^{(k)}(x) = \chi(rxr^{-1})$ for some $r$ such that $r^{-1}g_Kr = k$, and $\psi^{(l)}(x) = \psi(sxs^{-1})$ for some $s$ such that $s^{-1}g_Ls = l$.   

If $g_J$ is not in $KL$, then the sum is empty and $\langle \widehat{\chi} \otimes \widehat{\psi}, \widehat{\varphi} \rangle = 0$.
\end{Theorem}

\begin{Lemma}
Theorem \ref{actionfusion} and Theorem \ref{generalfusion} are equivalent.
\end{Lemma}

\begin{proof}
Note that $kl = g_J$ implies that $C(k) \cap C(l) \subset C(g_J)$. In \eqref{prettyformula}, the constant $\frac{|J|}{|G|}|C(k) \cap C(l)|$ is in fact equal to $[C(g_J) : C(k) \cap C(l)]^{-1}$, because $\frac{|J|}{|G|} = |C(g_J)|^{-1}$ by the orbit-stabilizer theorem. Moreover, there is an action of $C(g_J)$ on the set of pair $(k, l) \in KL$ such that $kl = g_J$: for $\gamma \in C(g_J)$, we can defined the action $\gamma \cdot (k, l) = (\gamma k \gamma^{-1}, \gamma l \gamma^{-1})$. This is well-defined since 
\[
(\gamma k \gamma^{-1})(\gamma l \gamma^{-1}) = \gamma kl \gamma^{-1} = \gamma g_J \gamma^{-1} = g_J.
\]
The stabilizer of $(k, l)$ under this action is $C(k) \cap C(l)$. Therefore, using the orbit-stabilizer theorem once again, we can remove the constant $[C(g_J) : C(k) \cap C(l)]^{-1}$ by only summing over the orbits. This is the same as summing over the set $O$ defined in the statement of Theorem \ref{actionfusion}.
\end{proof}

The formula \eqref{prettyformula} in Theorem \ref{generalfusion} can be written in another way:
\begin{equation}\label{fusionformula}
\langle \widehat{\chi} \otimes \widehat{\psi}, \widehat{\varphi} \rangle = \delta_{J \subset KL} \frac{|J|}{|G|} \sum_{\substack{(k, l) \in K \times L,\\kl = g_J}} \sum_{x \in C(k) \cap C(l)} \chi^{(k)}(x)\psi^{(l)}(x) \overline{\varphi(x)}.
\end{equation}
where we are summing over all $kl = g_J$ in the outer sum, and $\chi^{(k)}, \psi^{(l)}$ are defined as in Theorem \ref{generalfusion}. This is just expanding the inner product in (\ref{prettyformula}). Although it's more verbose, it's better suited to the applications.

We state the theorem in the special case where all centralizers are normal, which is proved in \cite{fusionrule}, since it's sometimes more convenient in applications.

\begin{Corollary}\label{fusion}
Denote $C(g_K) \cap C(g_L)= Q$. Then
\[
\langle \widehat{\chi} \otimes \widehat{\psi}, \widehat{\varphi} \rangle = \delta_{J \subset KL} \frac{|Q||J|}{|G|} \sum_{\substack{(k, l) \in K \times L, \\kl = g_J}} \langle \chi^{(k)} \psi^{(l)},  \varphi \rangle_Q.
\]
\end{Corollary}
\begin{proof}
If the centralizers are normal, then $C(k) = C(g_K)$ and $C(l) = C(g_L)$ for each $kl = g_J$, so the factor $|C(k) \cap C(l)|$ in \eqref{prettyformula} can be pulled out.
\end{proof}

\section{Fusion Rules via Mackey Theory}

In this section, we will prove Theorem \ref{generalfusion} in a slightly different form using Mackey theory. This alternative proof shows that the sum in Theorem \ref{generalfusion} is in fact a sum over certain double cosets. 

As before, let $J, K, L$ be conjugacy classes of $G$, $g_J, g_K, g_L$ be representatives of $J, K, L$ respectively, $H_1 = C(g_J), H_2 = C(g_K), H_3 = C(g_L)$ be the centralizers, and $(V_1, \pi_1), (V_2, \pi_2)$, $(V_3, \pi_3)$ be irreducible representations of $H_1, H_2, H_3$, respectively. 
Then the induced representations $V_i^G = \Ind_{H_i}^G V_i$ equipped with the $D(G)$-action as discussed in section \ref{DR} are irreducible representations of $D(G)$. The fusion coefficients are
\[
\dim \Hom_{D(G)}(V_2^G \otimes V_3^G, V_1^G).
\]
We would like to use Mackey theory to investigate $\Hom_{D(G)}(V_2^G \otimes V_3^G, V_1^G)$, but slightly differently than what we have done in section \ref{DR}.

We will see later that the use of Mackey theory involves understanding double cosets $H_3 \backslash G /H_2$. For notational convenience later on, we will choose a set of right coset representatives $\{\tau_k\}_{k \in K}$ for $H_2 \backslash G$ such that $\tau^{-1}_k g_K \tau_k = k$, and a set of left coset representatives $\{\sigma_l\}_{l \in L}$ for $G/H_3$ such that $\sigma_l g_L \sigma_l^{-1} = l$. Then, as $k, l$ vary through $K, L$, the product $\sigma_l^{-1}\tau_k^{-1}$ covers a complete set of representatives for double cosets in $H_3 \backslash G /H_2$.

\begin{Lemma}\label{sameaction}
For $(k, l), (k', l') \in K \times L$ such that $kl = k'l' = g_J$, $\sigma_l^{-1} \tau_k^{-1}$ represents the same double coset in $H_3 \backslash G / H_2$ as $\sigma_{l'}^{-1} \tau_{k'}^{-1}$ if and only if there is some $\gamma \in H_1$ such that $k' = \gamma k \gamma^{-1}$ and $l' = \gamma l \gamma^{-1}$.
\end{Lemma}
\begin{proof}
    The set of double cosets $H_3 \backslash G / H_2$ is in bijection with the orbits of $G/H_3 \times G/H_2$ under the action of $G$ by left multiplication. Under this bijection, $(\sigma_l H_3, \tau_k^{-1} H_2)$ corresponds to $H_3 \sigma_l^{-1}\tau_k^{-1}H_2$. Therefore, $\sigma_l^{-1} \tau_k^{-1}$ represents the same double coset in $H_3 \backslash G / H_2$ as $\sigma_{l'}^{-1} \tau_{k'}^{-1}$ if and only if there is some $\gamma \in G$ such that
    \[
    \gamma\sigma_l H_3 = \sigma_{l'}H_3 \quad \text{ and } \quad \gamma\tau_k H_2 = \tau_{k'}H_3
    \]
    We know that
    \[
    \gamma\sigma_{l} g_L \sigma_{l}^{-1}\gamma^{-1} = \gamma l \gamma^{-1}
    \]
    so $\gamma\sigma_{l}H_3 = \sigma_{\gamma l\gamma^{-1}}H_3$ and hence $\sigma_{\gamma l\gamma^{-1}}H_3 = \sigma_{l'}H_3$. The representatives $\sigma$'s are chosen beforehand, so we must have $l' = \gamma l\gamma^{-1}$. Similarly, $k' = \gamma k\gamma^{-1}$. We know
    \[
    g_J = k'l' = (\gamma l\gamma^{-1})(\gamma k\gamma^{-1}) = \gamma kl\gamma^{-1} = \gamma g_J \gamma^{-1}
    \]
    Hence $\gamma$ must be in $C(g_J) = H_1$. This is exactly what we wanted to prove.
\end{proof}

\begin{Theorem}\label{Mackeyfusion}
As above, let $H_1 = C(g_J)$, $H_2 = C(g_K)$, and $H_3 = C(g_L)$. Let $\tau_k$'s be a set of representative of $H_2 \backslash G$, and $\sigma_l$'s a set of representatives of $G/H_3$. $(V_i, \pi_i)$ is a irreducible representation of $H_i$, and $\chi_i$ is the character of $\pi_i$. Then
\[
\dim \Hom_{D(G)}(V_2^G \otimes V_3^G, V_1^G) =
\sum_{(k, l) \in O}\langle \chi_2^{(\tau_k)}\chi_3^{(\sigma_l^{-1})}, \chi_1 \rangle_{C(k) \cap C(l)}
\]
where $O$ is the set defined in Theorem \ref{actionfusion}.   
\end{Theorem}

\begin{proof}
Using Mackey's theorem on $\Res_{H_2}^G V_3^G$, we have that
\[
\Res_{H_2}^G V_3^G \cong \bigoplus_{s \in H_3 \backslash G / H_2} \Ind_{s^{-1}H_3s \cap H_2}^{H_2} \pi_3^{(s)} = \bigoplus_{s \in H_3 \backslash G / H_2} \left( \bigoplus_{r \in H_3 \backslash H_3 s H_2} r^{-1} \otimes V_3 \right)
\]
Induction commutes with direct sums, so we have that
\[
\Hom_{G}(V_2^G \otimes V_3^G, V_1^G) \cong \bigoplus_{s \in H_3 \backslash G / H_2} \Hom_{G}\left( \Ind_{H_2}^G \left(V_2 \otimes \Ind_{s^{-1}H_3s \cap H_2}^{H_2} \pi_3^{(s)} \right ),  V_1^G\right)
\]

We focus on a single summand. Applying Frobenius reciprocity, Mackey's theorem, and Frobenius reciprocity once again, we see that the summand for $s$ is 
\[
\bigoplus_{t \in H_2\backslash G/H_1} \Hom_{H_1 \cap t^{-1}H_2 t } \left( \pi_2^{(t)} \otimes \Ind_{t^{-1}s^{-1}H_3st \cap H_2}^{t^{-1}H_2 t} \pi_3^{(st)}, \pi_1 \right)
\]
As vector spaces, we have
\[
\pi_2^{(t)} \otimes \Ind_{t^{-1}s^{-1}H_3st \cap H_2}^{t^{-1}H_2 t} \pi_3^{(st)} = t^{-1} \otimes V_2 \otimes  \bigoplus_{r \in H_3 \backslash H_3 s H_2} t^{-1} r^{-1} \otimes V_3.
\]

For $t = \tau_k\rho_x$, $t^{-1}$ represents the same coset as $\tau^{-1}_{\rho_x^{-1}k\rho_x}$ in $G/H_2$ since $t^{-1}g_K t = \rho_x^{-1} \tau_k^{-1}g_K \tau_k \rho_x = \rho_x^{-1}k\rho_x$. Similarly, if we use the representative $r = \sigma_l^{-1}\tau_{k'}^{-1}$ then $t^{-1}r^{-1}$ represents the same coset as $\sigma_{yly^{-1}}$ in $G/H_3$ where $y = \rho_x^{-1}\tau_k^{-1}\tau_{k'}$. For any $g \in G$, the action of $\delta_g$ on $V_1$ is projecting to $\rho_g \otimes V_1$, so the action is $0$ unless $\rho_g$ represents the identity coset $H_1$. In other words, we only need to consider the action of $g = g_J$. This action on the source is given by
\[
\sum_{g_1 g_2 = g_J} \delta_{g_1}(t^{-1} \otimes V_2) \otimes \delta_{g_2}(t^{-1}r^{-1} \otimes V_3)
\]
and the action of $\delta_{g_i}$ is just projecting to the subspace $\tau^{-1}_{g_1} \otimes V_2$ and $\tau_{g_2} \otimes V_3$. Therefore, the action is zero unless 
\[
g_1 = \rho_x^{-1}k\rho_x \quad \text{ and } \quad g_2 = \rho_x^{-1}\tau_k^{-1}\tau_{k'}l (\rho_x^{-1}\tau_k^{-1}\tau_{k'})^{-1}.
\]
This is the same as saying that the action is non-zero only for $l$ such that 
\[
\rho_x^{-1}(k\tau_k^{-1}\tau_{k'}l\tau_{k'}^{-1}\tau_k)\rho_x = g_J
\]
Here both $k'$ and $k$ are fixed (we are varying only $r \in H_3 \backslash H_3 s H_2$, which is same as varying $l$), so we may rename $\tau_k^{-1}\tau_{k'}l\tau_{k'}^{-1}\tau_k$ as $l$. We see that the above condition is now simply $kl = x$. Clearly there is at most one $l$, i.e. at most one $r \in H_3 \backslash H_3 s H_2$ that can satisfy this condition. The action of $g_J$ then kills every summand other than
\[
\Hom_{H_1 \cap t^{-1}H_2t}(\rho_x^{-1}\tau_k^{-1} \otimes V_2 \otimes \rho_x^{-1}\sigma_l \otimes V_3, V_1) \text{ for } kl = x
\]

In summary, the $D(G)$-module morphisms in $\Hom_G(V_2^G \otimes V_3^G, V_1^G)$ is 
\[
\bigoplus_{\{(x, k, l) \in J \times K \times L \mid x = kl\}/\sim} \Hom_{H_1 \cap (\tau_k\rho_x)^{-1}H_2(\tau_k\rho_x)}(\rho_x^{-1}\tau_k^{-1} \otimes V_2 \otimes \rho_x^{-1}\sigma_l \otimes V_3, V_1)
\]
where $(x, k, l) \sim (x', k', l')$ if $\tau_k\rho_x$ represents the same double coset as $\tau_{k'}\rho_{x'}$, or $\sigma_{l}^{-1}\tau_k^{-1}$ represents the same double coset as $\sigma_{l'}^{-1}\tau_{k'}^{-1}$. 

Now we perform a re-indexing trick: re-name $k$ to be $\rho_x^{-1} k \rho_x$ and $l$ to be $\rho_x^{-1}l\rho_x$, then the new $\tau_k$ is the old $\tau_k\rho_x$, and the new $\sigma_l$ is the old $\rho_x^{-1}\sigma_l$. Notice that this re-indexing also eliminates the first condition: $\tau_{\rho_x^{-1} k \rho_x}\rho_{g_J}$ is just equal to $\tau_k\rho_x$. Therefore, we are left with
\begin{equation}\label{spacedecomp}
    \bigoplus_{\{(k, l) \in K \times L \mid g_J = kl\}/R} \Hom_{H_1 \cap \tau_k^{-1}H_2\tau_k}(\tau_k^{-1} \otimes V_2 \otimes \sigma_l \otimes V_3, V_1)
\end{equation}
where the relation $R$ is defined by $(k, l) R (k', l')$ if $\sigma_{l}^{-1}\tau_k^{-1}$ represents the same double coset as $\sigma_{l'}^{-1}\tau_{k'}^{-1}$. By Lemma \ref{sameaction}, we see that the indexing set is just the set $O$ in Theorem \ref{actionfusion}.

We have that $\tau_k^{-1}H_2\tau_k = C(k)$, and $H_1 \cap C(k)$ is equal to $C(l) \cap C(k)$ for $kl = g_J$. So now if we let $\chi_i$ be the character of $(V_i, \pi_i)$, then the dimension of \eqref{spacedecomp} is
\[
\sum_{(k, l) \in O} \langle \chi_2^{(\tau_k)}\chi_3^{(\sigma_l^{-1})}, \chi_1 \rangle_{C(k) \cap C(l)}.
\]
\end{proof}

Theorem \ref{Mackeyfusion} is just a re-statement of Theorem \ref{actionfusion}, but with a different point of view that focuses more on the structure of the induced representations. The proof is therefore of a much different flavor, highlighting that the summation is in fact over certain double coset representatives.

\section{Multiplicity Freeness for Dihedral Groups: $n$ odd}

In this and the next section, we compute fusion rules in the category of representations of $D(G)$ for $G$ Dihedral. It turns out these categories are multiplicity free. 

Notably, the Dihedral groups give interesting fusion rules where almost all simple objects appear as constituents. This behavior is very similar to the fusion rules of fusion algebras of Wess-Zumino-Witten conformal field theories, which we will discuss in the end of this paper.

The goal of this section to prove the following:
\begin{Theorem}\label{odd_mult_free}
Let $G$ be a Dihedral group $D_{2n}$ where $n$ is odd. Then the category of $D(G)$-modules is multiplicity free.
\end{Theorem}

We fix some representatives of conjugacy classes. Write $G = \langle x, y \mid x^n = y^2 = 1, yxy^{-1} = x^{-1}\rangle$, and $n = 2m+1$. Then the conjugacy classes of $G$ are $\{1\}$, $\{x^i, x^{-i}\}$ for $1 \leq i \leq m$, together with $\{yx^i \mid 0 \leq i \leq 2m\}$. We use the representatives $1$, $x^i$ for $1 \leq i \leq m$, and $y$. We see that the centralizer of $1$ is $G$, the centralizer of any of $x^i$ is the cyclic subgroup of order $n$ generated by $x$, and the centralizer of $y$ is $\{1, y\}$. We call simple objects parametrized by these three type of conjugacy classes type 1, type 2, and type 3, respectively. The proof of Theorem \ref{odd_mult_free} will be dealing with tensor products between different types of simple objects, so we do that one by one. 

\begin{Lemma}\label{t1t1}
The type-1-type-1 fusion rules are multiplicity free.
\end{Lemma}
\begin{proof}
A type-1 simple object is just a representation $V_{\pi}$ of $G = D_{2n}$ equipped with the $D(G)$ action. Thus, computing the fusion rule amounts to computing products of irreducible characters. For $n$ odd, we know that $D_{2n}$ has two one-dimensional characters: the trivial character $1$, and the sign character $s$ (whose value is $-1$ on the conjugacy class of $y$, and $1$ everywhere else.) There are $\frac{n-1}{2}$ two-dimensional irreducible characters $\chi_a$ for $a = 1, \cdots, \frac{n-1}{2}$. If $\zeta_n$ is a primitive $n$-th root of unity, let $\theta_a = \zeta_n^a + \zeta_n^{-a}$. The character $\chi_a$ sends $x^k$ to $\theta_{ka}$, and sends $y$ to $0$. We easily see that $s\chi_a = \chi_a$, so it remains to compute $\chi_a \chi_b$.

Notice that
\[\theta_a\theta_b = \theta_{a+b} + \theta_{a-b}.\]
Therefore, $\chi_a\chi_b$ sends $x$ to $\theta_{a+b} + \theta_{a-b}$, and sends $x^k$ to $\theta_{ka}\theta_{kb} = \theta_{k(a+b)} + \theta_{k(a-b)}$. So for $a \neq b$, we see that $\chi_a \chi_b = \chi_{a+b} + \chi_{a-b}$. For $a=b$, we have that $\chi_a^2$ sends $x^k$ to $\theta_{2ka} + 2$, and $y$ to $0$. If follows that $\chi_a^2$ contains a copy of $\chi_{2a}$, and the difference sends $y$ to $0$ and everything else to $2$. We then easily see that $\chi_a^2 = 1 + s + \chi_{2a}$. This covers all products, and they are indeed all multiplicity free.
\end{proof}

\begin{Lemma}\label{dihedraloddt1t2}
The type-1-type-2 fusion rules are all multiplicity free. 
\end{Lemma}
\begin{proof}
Here it is fine to use Corollary \ref{fusion}. We have $K = \{1\}$ and $L = \{x^i, x^{-i}\}$ for some $i$ fixed. Also, $C(g_K) = G = D_{2n}$, and $C(g_L) = C_n = \{1, x, \cdots, x^{n-1}\}$. Both of these are normal. We pick $g_L = x^i$. Let $\chi$ be some irreducible character of $C(g_K) = G$. Take $\psi$ to be a character of $C(g_L)$, which sends $x$ to some $n$-th root of unity. Let $J = L$. So $Q = G \cap C(g_L) = C(g_L) = \{1, x, \cdots, x^{n-1}\}$, and $KL = L$. So $J \subset KL$. We have $\frac{|Q||J|}{|G|} = \frac{n \cdot 2}{2n} = 1$. Let $\varphi$ be any character of $C(g_J)$. The $g_J$ we select is $x^i$, and there is only one way to write $x^i$ as a product of elements from $K$ and $L$: $x^i = 1 \cdot x^i$. Hence by Corollary \ref{fusion},
\[
\langle \widehat{\chi} \otimes \widehat{\psi}, \widehat{\varphi} \rangle = \langle \chi|_Q \psi, \varphi \rangle_Q.
\]
Now, if $\chi$ is one-dimensional, then its restriction to $Q = \{1, x, \cdots, x^{n-1}\}$ is always trivial, so in this case the inner product $\langle \chi|_Q \psi, \varphi \rangle_Q$ is $1$ if and only if $\varphi = \psi$. Now dimension considerations gives $V_{1, \chi} \otimes V_{x^i, \varphi} = V_{x^i, \varphi}$.

If $\chi$ is two dimensional, then it sends $x$ to some $\theta_a = \zeta_n^a + \zeta_n^{-a}$ for some primitive $n$-th root of unity. So its restriction to $Q$ does the same thing, i.e. it is a sum of of the character $x \mapsto \zeta_n^a$ and $x \mapsto \zeta_n^{-a}$ of $Q$. We know that $\psi$ is of the form $x \mapsto \zeta_n^b$ for some $b$, so 
\[\chi|_Q\psi = (x \mapsto \zeta_n^{a+b}) + (x \mapsto \zeta_n^{b-a})\]
Hence it contains one copy of $x \mapsto \zeta_n^{a+b}$ and one copy of $x \mapsto \zeta_n^{a+b}$. This shows that for $\pi$ two dimensional, $V_{1, \pi} \otimes V_{x^i, \theta}$ is a sum of one copy of $V_{x^i, \theta_{a+b}}$ and one copy of $V_{x^i, \theta_{b-a}}$. Note that $a + b$ is not equal to $b-a$ since $\chi$ being irreducible two dimensional implies $a \neq 0$. Hence type-1-type-2 fusion rules are all multiplicity free.
\end{proof}

\begin{Lemma}
The type-1-type-3 fusion rules are all multiplicity free.    
\end{Lemma}
\begin{proof}
This time we must use Theorem \ref{generalfusion}. Here, we have $K = \{1\}$ and $L = \{yx^j \mid j=0, \cdots, n-1\}$. Also, $C(g_K) = G = D_{2n}$, and $C(g_L) = C_2 = \{1, y\}$. As before we pick $g_L = y$. Let $\chi$ be some irreducible character of $C(g_K) = G$. Take $\psi$ to be a character of $C(g_L)$, which can only be one-dimensional. Let $J = L$, so $J \subseteq KL = L$. Thus $\delta_{J \subset KL} = 1$, and $\frac{|J|}{|G|} = \frac{n}{2n} = \frac{1}{2}$. Let $\varphi$ be any character of $C(g_J)$, which again can only be one-dimensional. The $g_J$ we select is just $y$, and there is only one way to write $y$ as a product of elements from $K$ and $L$: $y = 1 \cdot y$. So we must have $k = 1, l = y, r = 1$, and $s = 1$. Then in \eqref{fusionformula}, we are only summing over elements in $C(g_K) \cap C(g_L) \cap C(g_L) = \{1, y\}$. So we have
\[
\langle \widehat{\chi} \otimes \widehat{\psi}, \widehat{\varphi} \rangle = \frac{1}{2} (\chi(1)\psi(1)\varphi(1) + \chi(y)\psi(y)\varphi(y))
\]
If $\chi$ is a one-dimensional character, then 
\[
\chi|_Q\psi = \begin{cases}
1, & \chi|_Q = \psi \\
s, & \chi|_Q \neq \psi
\end{cases}
\]
In either case, the product $\chi(y)\psi(y)\varphi(y)$ is $1$ if $\varphi(y) = \chi(y)\psi(y)$ and $-1$ if not (there are only two choices of $\varphi(y)$). So this shows that the tensor product of a type-1 simple object $V_{1, \chi}$ and a type-3 simple object $V_{y, \psi}$ contains only one of $V_{y, 1}$ and $V_{y, s}$. By considering dimensions we see that this is the only irreducible summand.

Now suppose $\chi$ is a two dimensional irreducible character, which necessarily sends $y$ to $0$ and $1$ to $2$. Thus $\chi(y)\psi(y)\varphi(y)$ is $0$, and $\chi(1)\psi(1)\varphi(1) = 2$, so $\langle \widehat{\chi} \otimes \widehat{\psi}, \widehat{\varphi} \rangle = 1$ no matter which character $\varphi$ is. This means that $V_{1, \chi} \otimes V_{y, \psi}$ contains one copy of $V_{y, 1}$ and one copy of $V_{y, s}$. Again by dimension consideration, these are all the irreducible summands. Hence type-1-type-3 fusion rules are all multiplicity free.
\end{proof}

\begin{Lemma}
The type-2-type-2 fusion rules are all multiplicity free.  
\end{Lemma}
\begin{proof}
We deal with simple objects of the form $V_{x^i, \pi}$. For any $i$, $C(x^i)$ is the cyclic group of order $n$, which we denote by $C_n$. The group $C_n$ has $n$ irreducible representations, each mapping the generator to some $n$-th root of unity. Let them be $\rho_1, \cdots, \rho_n$. Using the same notation as above, we have $K = L = \{x^i, x^{-i}\}$, and $C(g_K) = C(g_L) = \{1, x, \cdots, x^{n-1}\}$ is normal in $D_{2n}$. Hence, the formula in Theorem \ref{fusion} works. Letting $J = \{x^i, x^{-i}\}, K = \{x^j, x^{-j}\}$ and $J = \{x^{i+j}, x^{-i-j}\}$, we easily see that the tensor product $V_{x^i, \rho_k} \otimes V_{x^j, \rho_l}$ contains a copy of $V_{x^{i+j}, \rho_{k+l}}$. Moreover, if $i \neq j$, then the same computation shows it also contains a copy of $V_{x^{i-j}, \rho_{k-l}}$. 

If $i = j$, the $i+j$ computation remains the same since it is not possible that $i+j = 2i = n$. This gives an irreducible consituent $V_{x^{2i}, \rho_{k+l}}$ of dimension $2$. For the $i-j$ computation, we let $J = \{1\}$. Let $\chi, \psi$ be two one-dimensional characters of $C(g_K)$, and let them be $\rho_k$ and $\rho_l$. Then $Q = \{1, x, \cdots, x^{n-1}\}$, and $\frac{|Q||J|}{|G|} = \frac{1}{2}$. There are two ways to write $1$: $1 = x^i x^{-i} = x^{-i}x^i$. In the first decomposition, we twist the second character $\psi = \rho_l$ by $y$, which has the effect of making it $\rho_{-l}$. So $\chi\psi = \rho_k\rho_{-l} = \rho_{k-l}$, which sends $x$ to $\zeta_n^{k-l}$. If $k \neq l$, then the two dimensional irreducible character $\varphi = (x \mapsto \theta_{k-l})$ of $G$ restricts to $Q$ is $\rho_{k-l} + \rho_{l-k}$, so $\langle \chi\psi^{(y)}, \varphi \rangle_Q = 1$. For the other way of decomposition, we twist $\chi$ by $y$, which makes it $\rho_{-k}$. So again $\langle \chi^{(y)}\psi, \varphi \rangle_Q = 1$. Hence the multiplicity of $\varphi = (x \mapsto \theta_{k-l})$ in $V_{x^i, \rho_k} \otimes V_{x^i, \rho_l}$ is $\frac{1}{2}(1+1) = 1$.

If $k = l$, then $\rho_{k-l} = \rho_{l-k} = 1$, and we know that both one dimensional irreducible characters of $G$ restrictts to $1$. Thus $V_{x^i, \rho_k} \otimes V_{x^i, \rho_k}$ contains $V_{1, 1} \oplus V_{1, s}$. So in particular a formula for squaring type 2 simple objects is
\[V_{x^i, \rho_k}^{\otimes 2} = V_{1, 1} \oplus V_{1, s} \oplus V_{x^{2i}, \rho_{2k}}.\]
In any case we see that the $i-j$ computation still gives irreducible constituents of dimension $2$ in total. Thus we found all irreducible constituents in all cases, and they are all multiplicity free.
\end{proof}

\begin{Lemma}
The type-2-type-3 fusion rules are all multiplicity free.    
\end{Lemma}
\begin{proof}
We use Theorem \ref{generalfusion} in the same way as before. Here, we have $K = \{x^i, x^{-i}\}$ for some fixed $i$ and $L = \{yx^j \mid j=0, \cdots, n-1\}$. Let $J = L$. Here $C(g_K)$ is indeed normal, so elements in the inner sum of \eqref{fusionformula} must be in $C(g_K) \cap C(g_L) = 1$, and $KL = L$. Thus $J \subset KL$, and $\frac{|Q||J|}{|G|} = \frac{1}{2}$. The $g_J$ we select is just $y$, and to write $y$ as a product of elements from $K$ and $L$, notice that $yx^j = x^{-j}y$. So we have the following ways:
\[
y = x^i yx^i = x^{-i} yx^{-i}
\]
So the sum has two terms. Since we have seen that the inner sum can only be at $1$, and $\chi, \psi, \varphi$ are all one-dimensional, the summand can only be $1$. Thus we see that 
\[
\langle \widehat{\chi} \otimes \widehat{\psi}, \widehat{\varphi} \rangle = \frac{1}{2} \times 2 = 1
\]
for any choice of $\chi$, $\psi$ and $\varphi$. In other words, for any simple object $V_{x^i, \pi}$ of type 2 and $V_{y, \theta}$ of type 3, their tensor product contains one copy of $V_{y, 1}$ and one copy of $V_{y, s}$ where $s$ is the representation of $\{1, y\}$ that sends $y$ to $-1$.

On the other hand, $V_{x^i, \pi}$ has dimension $2$, and $V_{y, \theta}$ has dimension $n$, so their tensor product has dimension $2n$. Since $V_{y, 1} \oplus V_{y, s}$ already has dimension $2n$, they are all the irreducible summands in $V_{x^i, \pi} \otimes V_{y, \theta}$. This shows that the type-2-type-3 fusion rules are all multiplicity free.
\end{proof}

\begin{Lemma}\label{oddt3t3}
The type-3-type-3 fusion rules are all multiplicity free.  
\end{Lemma}
\begin{proof}
In the notation of Theorem \ref{generalfusion}, we let $K = L = \overline{y} = \{yx^j \mid j = 0, \cdots n-1\}$, and $C(g_K) = C(g_L) = \{1, y\}$. Let $\chi, \psi$ be irreducible characters of $\{1, y\}$. They can only be the trivial representation or the sign representation.

First let $J = \{x^i, x^{-i}\}$ for some $1 \leq i \leq m$. Then $KL$ contains $J$ since $x^i = y(yx^i)$. Also $C(g_J) = \{1, x, \cdots, x^{n-1}\}$. We see that the inner sum of the formula in Theorem \ref{generalfusion} has only one term (i.e. at $1$). We pick $g_J = x^i$. Let $\varphi$ be any irreducible character of $C(g_J)$, which sends $x$ to some $n$-th root of unity. We have $\frac{|Q||J|}{|G|} = \frac{1 \cdot 2}{2n} = \frac{1}{n}$. We know that $yx^jyx^k = x^{k-j}$, so there are $n$ ways to write $g_J = x^i$ as a product of elements from $K$ and $L$:
\[x^i = (yx^j)(yx^i) \text{ for } k-j \equiv i \text{ mod } n.\]
The term $\langle \chi^{(r)} \psi^{(s)},  \varphi \rangle_Q$ is always $1$ since $Q$ is the trivial group and all these characters are one-dimensional. Thus we see that
\[
\langle \widehat{\chi} \otimes \widehat{\psi}, \widehat{\varphi} \rangle = \frac{1}{n} \cdot n = 1
\]
for each $1 \leq i \leq m$.

Next, let $J = \{1\}$, and $g_J = 1$. Certainly $KL$ contains $J = \{1\}$. We have $\frac{|J|}{|G|} = \frac{1}{2n}$. Let $\varphi$ be some two-dimensional irreducible character of $C(g_J) = G = D_{2n}$. There are $n$ ways to write $1$ as a product of elements from $K$ and $L$, since each $yx^j$ is its own inverse. Now, referring to \eqref{fusionformula}, in the sum we always have $k = l = yx^i$ for some $0 \leq i \leq n-1$, and $r = s = x^{a}$ where $2a \equiv i$ mod $n$ (because $x^{-a}yx^{a} = yx^{2a}$). For such $r$, we see that $rC(k)r^{-1} = x^a\{1, y\}x^{-a} = \{1, yx^{-i}\}$. Any two dimensional character $\varphi$ of $G$ sends $1$ to $2$, and sends $y$ to $0$, so each term is equal to $1 \times 2 + 0 = 2$. We have $n$ terms, so the result is $\frac{1}{2n} \cdot n \cdot 2 = 1$. This means that the tensor product of two type-3 simple objects contains one copy of each two-dimensional type-1 simple object.

If $\varphi$ is a one-dimensional character of $G$, then it is either the trivial character or the sign character. We see that the sum contains $n$ terms equaling $1$ as before, but now $\chi^{(r)}(yx^{-i}) = \chi(x^{a}yx^{-i}x^{-a}) = \chi(yx^{-2a-i}) = \chi(y)$, and similarly $\psi^{(r)}(yx^{-i}) = \psi(y)$, and $\varphi^{(r)}(yx^{-i}) = \varphi(y)$, so the result depend on whether $\chi, \psi, \varphi$ send $y$ to $1$ or $-1$. When $\chi = \psi$, we see that we get $n$ for $\varphi = 1$, so in total we have $\frac{1}{2n}(n+n) = 1$; we get $-n$ for $\varphi = s$ sending $y$ to $-1$, so in total the sum is $n-n = 0$. Therefore we see that for $\chi = \psi$, the tensor product contains one copy of $V_{1, 1}$ (and $0$ copies of the one with the sign representation). When $\chi \neq \psi$, we easily see that the results are opposite.

In summary, the tensor product of two type-3 simple objects contains one copy of each simple object parametrized by some $x^i$ (giving dimension $2mn$ in total), and one copy of each simple object parametrized by $1$ and some 2-dimensional irreducible representation (giving dimension $2m$), and one copy of the one dimensional simple object as discussed above. This accounts for all simple summands, since the dimension of the tensor product is $n^2$, and we accounted for $2mn + 2m + 1=n^2$.
\end{proof}

This finishes the proof of Theorem \ref{odd_mult_free}. Notice that the type-3-type-3 fusion rules give a sum involving most of the simple objects.

\section{Multiplicty Freeness for Dihedral Groups: $n$ even}\label{diheven}
Next, we deal with the case where $n = 2m$ is even.

Again, we fix the notations. For $G = D_{2n} = \langle x, y \mid x^n = y^2 = 1, yxy^{-1} = x^{-1} \rangle$ with $n$ even, the conjugacy classes are $\{1\}$, $\{x^m\}$, $\{x^i, x^{-i}\}$ for $1 \leq i < m$, $\{yx^j \mid 0 \leq j \leq n-2, j \text{ even}\}$ and $\{yx^j \mid 1 \leq j \leq n-1, j \text{ odd}\}$. They have size $1, 1, 2, m, m$, respectively. Indeed, we have $2n = 1 + 1 + 2(m-1) +m + m$. We fix the representatives to be $1, x^m, x^i, y, yx$, respectively. Moreover, using the relation $yx^iy^{-1} = x^{-i}$, we have that $C(1) = G$, $C(x^m) = G$, $C(x^i) = \{1, \cdots, x^{n-1}\}$, $C(y) = \{1, y, x^m, yx^m\} \cong C_2 \times C_2$, and $C(yx) = \{1, yx, x^m, yx^{m+1}\} \cong C_2 \times C_2$. We call simple objects parametrized by these conjugacy classes type $1$, type $1'$, type $2$, type $3$, and type $3'$.

To start, we discuss the irreducible characters of $G = D_{2n}$ with $n$ even. The commutator subgroup of $G$ is $\langle x^2 \rangle$, so the abelianization $G/[G, G]$ has order $4$, with representatives $\{1, y, x, yx\}$, which is isomorphic to $\{1, y\} \times \{1, x\} \cong C_2 \times C_2$. Therefore, there are $4$ one-dimensional characters of $G$ pulled back from the characters of $C_2 \times C_2$. We call them $(1, 1), (1, s), (s, 1)$ and $(s, s)$, where $s$ is the sign character as usual, the first slot is a character of $\{1, y\}$ and the second is a character of $\{1, x\}$. So $(s, 1)$ sends $y$ to $-1$ and $x$ to $1$, etc. 

Similar to the odd case, there are two dimensional characters $\chi_a$ sending $y$ to $0$, $yx$ to $0$, and $x^k$ to $\theta_{ka} = \zeta_n^{ka} + \zeta_n^{-ka}$ for $1 \leq a \leq m - 1$. We note that $\chi_m$ is not irreducible, as it decomposes as $(1, s) + (s, s)$. These are all irreducible characters.

\begin{Lemma}\label{event1t1}
The type-1-type-1 fusion rules are multiplicity free.
\end{Lemma}
\begin{proof}
Again it suffices to compute the product of characters of $D_{2n}$. For any two-dimensional irreducible character $\chi_a$, we see that $(s, 1)\chi_a = \chi_a$ since $\chi_a$ sends $y$ and $yx$ to $0$. We have that $(1, s)\chi_a$ sends $x$ to $-\theta_a$, and for $n$ even, $-1 = \zeta_n^{m} = \zeta_n^{-m}$. So $-\theta_a = \zeta_n^{-m+a} + \zeta_n^{m-a} = \theta_{m - a}$. It follows that $(1, s)\chi_a = \chi_{m-a}$. Likewise $(s, s)\chi_a = \chi_{m-a}$ since $\chi_a(y) = \chi_a(yx) = 0$.

Next, if $\chi_a$ and $\chi_b$ are irreducible two dimensional characters, then $\chi_a(x)\chi_b(x) = \theta_{a+b} + \theta_{a-b}$, as computed in Lemma \ref{t1t1}. If $a \neq b$, then $\chi_a \chi_b = \chi_{a+b} + \chi_{a - b}$. It is possible that $a + b = m$, in which case $\chi_a\chi_b = (1, s) + (s, s) + \chi_{a-b}$. If $a = b$, then $\chi_a^2$ sends $x^k$ to $\theta_{2ka} + \theta_0 = \theta_{2ka} + 2$. Thus $\chi_a^2$ contains $1 + (s, 1)$ and $\chi_{2a}$ as irreducible constituents. If $2a \neq m$, then $\chi_{2a}$ is irreducible, and if $2a = m$, then $\chi_{2a} = \chi_m$ decomposes as $(1, s) + (s, s)$. In summary, $\chi_{a}^2 = 1 + (s, 1) + \chi_{2a}$ if $2a \neq m$, and $\chi_{a}^2 = 1 + (s, 1) + (1, s) + (s, s)$ if $2a = m$. This shows all such fusion rules are multiplicity free.
\end{proof}

We may use the same computation for type-$1$-type-$1'$ and type-$1'$-type-$1'$ fusion rule, but there is a subtlety: type 1 and type 1' simple objects are isomorphic as $\C[G]$-modules, but not isomorphic as $D(G)$-module. One can check that type-$1$-type-$1'$ tensor products result in a sum of type $1'$ objects, and type-$1'$-type-$1'$ tensor products result in a sum of type 1 objects (the summands depend on the irreducible characters as in Lemma \ref{event1t1}). Consequently, we have

\begin{Corollary}
The type-$1$-type-$1'$ and type-$1'$-type-$1'$ fusion rules are all multiplicity free.
\end{Corollary}

\begin{Lemma}
The type-1-type-2 fusion rules are multiplicity free.
\end{Lemma}
\begin{proof}
The proof is similar to that of \ref{dihedraloddt1t2}. We will just state the fusion rules here.
Let $\rho_a: x \mapsto \zeta_n^a$ for some $a$, where $\zeta_n$ is a primitive $n$-th root of unity. Then we have
\[
V_{1, (1, 1)} \otimes V_{x^i, \rho_b} = V_{1, (s, 1)} \otimes V_{x^i, \rho_b}= V_{x^i, \rho_b} \quad \text{and} \quad
V_{1, (1, s)} \otimes V_{x^i, \rho_b} = V_{1, (s, s)} \otimes V_{x^i, \rho_b}= V_{x^i, \rho_{m+b}}
\]

When $\chi$ is a two-dimensional irreducible character of $G$, we have
\[
V_{1, \chi_a} \otimes V_{x^i, \rho_b} = V_{x^i, \rho_{a+b}} \oplus V_{x^i, \rho_{a-b}}
\]
which is also multiplicity free.
\end{proof}

\begin{Lemma}\label{dihedralevent2t2}
The type-2-type-2 fusion rules are multiplicity free.
\end{Lemma}
\begin{proof}
Let $K = \{x^i, x^{-i}\}$ and $L = \{x^j, x^{-j}\}$ for some $i, j \neq m$. Suppose $i \neq j$ and $x^{i + j} \neq x^m$. First let $J = \{x^{i+j}, x^{-i-j}\}$. In view of Corollary \ref{fusion}, we have $Q = C_n = \{1, x, \cdots, x^{n-1}\}$. Let $\chi, \psi, \varphi$ be characters of $C(g_K) = C(g_L) = C(g_J) = Q$. Then we have $\frac{|Q||J|}{|G|} = \frac{n \cdot 2}{2n} = 1$, and there is only one way to write $x^{i+j}$ as a product of two elements from $K$ and $L$, namely $x^{i + j} = x^ix^j$. Thus the sum contains only one term. Therefore, we see that the term is $1$ if $\varphi = \chi\psi$ and $0$ otherwise. This means $V_{x^i, \chi} \otimes V_{x^j, \psi}$ contains a copy of $V_{x^{i+j}, \chi\psi}$. Similarly, if we consider $J = \{x^{i-j}, x^{j-i}\}$ (using the assumption that $i \neq j$), we will see that $V_{x^i, \chi} \otimes V_{x^j, \psi}$ contains a copy of $V_{x^{i-j}, \chi\psi}$. Note that $V_{x^i, \chi} \otimes V_{x^j, \psi}$ is has dimension $4$, and we found two irreducible constituents of dimension $2$, so we found all irreducible constituents.

Now if $i = j$, the set $J = \{x^{i-j}, x^{j-i}\}$ is just $J = \{1\}$. In this case, $Q = C(g_K) \cap G = C(g_K) = \{1, \cdots, x^{n-1}\}$ as before, but now $\frac{|Q||J|}{|G|} = \frac{1}{2}$. Moreover, there are two ways to write $1$, namely $x^ix^{-j}$ and $x^{-i}x^j$. The first twists the character $\chi$ by the identity and the character $\psi$ by $y$, and the second twists the character $\chi$ by $y$ and the character $\psi$ by identity. If $\chi^{(y)}\psi$ is equal to the trivial character, then since $\chi\psi^{(y)}$ sends $x^k$ to $(\chi\psi)(x^{-1})$, it is also the trivial character. So the sum for $\varphi$ such that $\varphi|_Q = 1$ is $\frac{1}{2}(1+1) = 1$. There are two irreducible $\varphi$ such that $\varphi|_Q$ is $1$ (i.e. $(1, 1)$ and $(s, 1)$), so we have that $V_{x^i, \chi} \otimes V_{x^j, \psi}$ contains $V_{1, (1, 1)} \oplus V_{1, (s, 1)}$. If $\chi^{(y)}\psi$ is equal to the sign characters, then similarly $V_{x^i, \chi} \otimes V_{x^j, \psi}$ contains $V_{1, (1, s)} \oplus V_{1, (s, s)}$. 

If $\chi^{(y)}\psi$ is not equal to either the trivial character or the sign character, then say $\chi^{(y)} = \rho_a$ and $\psi = \rho_b$. Then $\chi^{(y)}\psi = \rho_{a+b}$. Let $\varphi$ be the two dimensional character of $G$ sending $x^k$ to $\theta_{k(a+b)}$, which decomposes as $\rho_{a+b}$ and $\rho_{a-b}$ after restricting to $Q$. Then we see that $\langle \chi^{(y)}\psi, \varphi \rangle = 1$. On the other hand, $\chi\psi^{(y)}$ sends $x$ to $\zeta_a\zeta_{-b}$, so it is $\rho_{a-b}$. Hence $\langle \chi\psi^{(y)}, \varphi \rangle = 1$ too. In summary, for this $\varphi$, its multiplicity in the tensor product is $\frac{1}{2}(1 + 1) = 1$. In any case, we see that when $i = j$, the fusion is still multiplicity free.

An identical computation can be carried out for the situation where $j + i = m$, and of course $i = j$ and $i + j = m$ can happen at the same time. But we always have a multiplicity free fusion rule since the above computed irreducible constituents cover the correct dimension.
\end{proof}

We are now left with fusion rules involving type 3 and type $3'$ simple objects. We will only consider type 3 simple objects, since the ones for type $3'$ will be mostly identical. Note that type 3 (and type $3'$) simple objects have dimension $\frac{n}{2} = m$.

\begin{Lemma}\label{event1t3}
The type-1-type-3 fusion rules are multiplicity free.
\end{Lemma}
\begin{proof}
Let $K = \{1\}$ and $L = \{yx^j \mid j \text{ even}\}$. The representative for $L$ is chosen to be $y$. Then $C(g_L) = \{1, y, x^m, yx^m\}$. Let $\chi$ be a one dimensional character of $G$, and $\psi$ be some character of $C(g_L)$. We denote the characters of $C(g_L) = C_2 \times C_2$ by $(1, 1), (t, 1), (1, t), (t, t)$, where $(1, t)$ sends $x^m$ to $-1$, etc. Let $J = L$. Here we must use Theorem \ref{generalfusion} as $C(g_L)$ is not normal. Clearly $J \subset KL$, and $\frac{|J|}{|G|} = \frac{|m|}{|2n|} = \frac{1}{4}$. There is only one way to write $y$ as a product of two elements in $K$ and $L$, namely $y = 1 \cdot y$. So all twists are by the identity, so in \eqref{prettyformula} we are summing over $x \in C(g_K) \cap C(g_L) \cap C(g_J) = C(g_L) = \{1, y, x^m, yx^m\}$. In other words, we are computing
\[
\frac{1}{4} \cdot 4 \cdot \langle \chi\psi, \varphi \rangle_{C(g_L)} = \langle \chi\psi, \varphi \rangle_{C(g_L)}
\]
Recall that we are assuming $\chi$ is one dimensional, so it's restriction to $C(g_L)$ is irreducible. This means $V_{1, \chi} \otimes V_{y, \psi}$ contains one copy of $V_{y, \chi|_{C(y)}\psi}$. Both of these are of dimension $4$, so this is the fusion rule, which is multiplicity free.

Now suppose $\chi$ is a two dimensional character $\chi_a$ sending $x^k$ to $\theta_{ka}$. Then the restriction of $\chi$ to $C(g_L)$ sends $x^m$ to $\theta_{km} = \zeta_{n}^{am} + \zeta_{n}^{-am} = (-1)^a + (-1)^{-a} = \pm 2$. Of course it sends $y, yx^m$ to $0$ too, and $1$ to $2$. Thus it is $(1, 1) + (t, 1)$ or $(1, t) + (t, t)$. This means 
\[
V_{1, \chi} \otimes V_{y, \psi} = V_{y, (1, 1)} \oplus V_{y, (t, 1)} \text{ or } V_{y, (1, t)} \oplus V_{y, (t, t)}
\]
The tensor product is of dimension $8$, which equals to the dimension of the sum. Hence this is the complete fusion rule, and it is multiplicity free.
\end{proof}

Note that the character theoretic computations for type $1'$ and type $3'$ are the same, so we will not do them again.

\begin{Lemma}\label{dihedralevent2t3}
The type-2-type-3 fusion rules are multiplicity free.
\end{Lemma}
\begin{proof}
Let $K = \{x^i, x^{-i}\}$ for some $i \neq m$. Let $L = \{yx^j \mid j \text{ even}\}$ as above. Suppose first $i$ is odd. Let $J = \{yx^j \mid j \text{ odd}\}$. Recall that $C(g_J) = \{1, yx, x^m, yx^{m+1}\}$. We have $\frac{|J|}{|G|} = \frac{1}{4}$. Then $yx^{-i}$ is in side $J$. We may pick the representative of $J$ to be $yx^i$ here for convenience. Then there are two ways to write $yx^i$, namely
\[
yx^{-i} = x^{i}y = x^{-i}(yx^{-2i}).
\]
Corresponding to the first way of decomposition, in \ref{generalfusion} we have summing over all elements $x$ in $C(g_K) \cap C(g_L) \cap C(g_J) = \{1, x^m\}$. A character $\chi$ of $C(g_K)$ is of the form $\rho_a$ for some $a$, so it sends $1$ to $1$ and $x^m$ to $-1$, and characters $\psi, \varphi$ of $C(g_L)$ and $C(g_J)$ are all one dimensional, being products of the trivial and the sign character. Thus
\[
\chi(1)\psi(1)\varphi(1) + \chi(x^m)\psi(x^m)\varphi(x^m) = 1 - \psi(x^m)\varphi(x^m)
\]
No matter what $\psi$ is, there are two choices of $\varphi$ such that $\psi(x^m)\varphi(x^m) = -1$. Call them $\varphi_1$ and $\varphi_2$. For these two characters, we compute the other sum given by $yx^{-i} = x^{-i}(yx^{-2i})$. Here, $x^{-i} = yx^{i}y^{-1}$, and $yx^{-2i} = x^iyx^{-i}$. So we are summing over elements in 
\[
y^{-1}C(g_K)y \cap x^{-i}C(g_L)x^i \cap C(g_J) = \{1, x^m\}.\]
So the result is in fact the same. We conclude that for $\varphi = \varphi_1, \varphi_2$, their multiplicities are both $\frac{1}{4}(1+1+1+1) = 1$. Hence for $i$ odd,
\[
V_{x^i, \chi} \otimes V_{y, \psi} = V_{yx, \varphi_1} \otimes V_{yx, \varphi_2}.
\]
The equality is again by dimension considerations. This is multiplicity free.

Now suppose $i$ is even. This time let $J = L$. We conveniently choose the representative $yx^{-i}$ again. The computation will be exactly the same as above, since the intersection in both decompositions are still $\{1, x^m\}$. This means for $i$ even,
\[
V_{x^i, \chi} \otimes V_{y, \psi} = V_{y, \varphi_1} \otimes V_{y, \varphi_2}.
\]
where $\varphi_1, \varphi_2$ are characters of $\{1, y, x^m, yx^m\}$ such that $\psi(x^m)\varphi_1(x^m) = \psi(x^m)\varphi_2(x^m) = -1$.
\end{proof}

Finally, we compute the fusion rules between type 3 simple objects, which are the most interesting fusion rules.

\begin{Lemma}\label{event3t3}
The type-3-type-3 fusion rules are multiplicity free.
\end{Lemma}
\begin{proof}
The proof follows a similar procedure as Lemma \ref{oddt3t3}. We list the results here. Let $\chi, \psi$ be irreducible characters of the centralizer $\{1, y, x^m, yx^m\}$. When $m$ is odd, there is 1 type-1 one-dimensional simple constituent, $\frac{m-1}{2}$ type-1 two-dimensional simple constituents, $\frac{m-1}{2}$ type-2 simple constituents. 

When $m$ is even, if $\chi\psi$ is $(1, 1)$ or $(t, 1)$ then there are $2$ type-1 one-dimensional simple constituents, $\frac{m}{2} - 1$ type-1 two-dimensional simple constituents, $m$ type-$1'$ one-dimensional simple constituents, and $\frac{m}{2} - 1$ type-2 simple constituents. If $\chi\psi$ is $(1, t)$ or $(t, t)$, then there are $\frac{m}{2}$ type-1 two-dimensional simple constituents, $m$ type-$1'$ one-dimensional simple constituents, and $\frac{m}{2} - 1$ type-2 simple constituents.
\end{proof}

Note that the type-3-type-3 fusion rules give a sum that contains $o(m^2)$ summands. As such, it is a rich fusion rule.

Combining all the above lemmas, we have proved
\begin{Theorem}
Let $G$ be a Dihedral group $D_{2n}$ where $n$ is even. Then the category of $D(G)$-modules is multiplicity free.
\end{Theorem}

\section{Dicyclic Groups: $n$ even}

We now begin the study of another family of groups, the Dicyclic groups. In this and the next section we will prove that the category of finite dimensional representations of $D(Q_{4n})$ has multiplicity free fusion rules.

The dicyclic group $Q_{4n}$ of degree $n$ has the presentation
\[
\langle x, y \mid x^{2n} = 1, y^2 = x^n, yxy^{-1} = x^{-1} \rangle
\]

The conjugacy classes of $Q_{4n}$ are:
\begin{enumerate}
    \item Type 1, $1'$: $\{1\}, \{x^n\}$;
    \item Type 2: $\{x^i, x^{-i}\}$ for $i = 1, \cdots, n-1$;
    \item Type 3, $3'$: $\{yx^j \mid j \text{ even}\}$, $\{yx^j \mid j \text{ odd}\}$.
\end{enumerate}
As usual, we need to know the centralizers of each conjugacy class. For type 1 and $1'$, the centralizer is all of $Q_{4n}$. For type 2, the centralizer is the subgroup generated by $x$, i.e. $\{1, x, \cdots, x^{2n-1}\}$. Finally, the centralizer of $y$ is $\{1, y, x^n = y^2, yx^n = y^3\}$, and the centralizer of $yx$ is $\{1, yx, x^n = y^2, yx^{n+1} = y^3x\}$. Note that all the above hold true without assuming $n$ is even, so we will refer to them in the next section as well.

It is a fact that when $n$ is even, the character table of $Q_{4n}$ is the same as the character table of $D_{4n}$. Moreover, the centralizers of each conjugacy is identical to what we saw for $D_{4n}$, so in fact the computations in all lemmas in Section \ref{diheven} remain true here. In other words, they have the same fusion rules:
\begin{Theorem}
    The Grothendieck rings of $\mathrm{Rep}_f {D(Q_{4n})}$ and $\mathrm{Rep}_f {D(D_{4n})}$ are isomorphic.
\end{Theorem}
We then easily obtain that
\begin{Corollary}
Let $G$ be a Dicyclic group $Q_{4n}$ where $n$ is even. Then the category of $D(G)$-modules is multiplicity free.    
\end{Corollary}

We remark here that the generalized quaternion groups are dicyclic groups $Q_{4n}$ with $n$ a power of $2$. In particular, the usual quaternion group $Q_8$ is one example. So for all generalized quaternion groups, their Drinfeld doubles have multiplicity free representation categories.

\section{Dicyclic Groups: $n$ odd}

In this section we deal with $Q_{4n}$ where $n$ is odd. Many of the proofs will be omitted since they are similar to those in the previous sections.

There are $n-1$ irreducible two characters $\chi_a$ for $1 \leq a \leq n-1$ and $4$ one-dimensional characters. We name the four one-dimension characters by $1, t_i, t_{-1}, t_{-i}$ respective, the subscript indicating what $y$ is sent to.

Now we begin the computation of fusion rules.

\begin{Lemma}%\label{event1t1}
The type-1-type-1 fusion rules are multiplicity free.
\end{Lemma}
\begin{proof}
Similar to the Dihedral cases, it suffices to compute the product of characters of $Q_{4n}$. The computation is trivial if one of the character is one-dimensional .For two-dimensional characters $\chi_a, \chi_b$, it is easy to obtain the following results: 

If $a \neq b$, then $\chi_a \chi_b = \chi_{a+b} + \chi_{a - b}$. It is possible that $a + b = n$, in which case $\chi_a\chi_b = t_{i} + t_{-i} + \chi_{a-b}$. In this case $a \neq b$ since $n$ is odd. If $a = b$, then $\chi_a^2$ sends $x^k$ to $\theta_{2ka} + \theta_0 = \theta_{2ka} + 2$. Thus $\chi_a^2$ contains $1 + t_{-1}$ and $\chi_{2a}$ as irreducible constituents. It is not possible for $2a$ to be equal to $n$ since $n$ is odd. This shows all such fusion rules are multiplicity free.
\end{proof}

\begin{Lemma}
The type-1-type-2 fusion rules are multiplicity free.
\end{Lemma}
\begin{proof}
The proof is similar to that of Lemma \ref{dihedraloddt1t2}. We will just state the fusion rule here:  

Let $\rho_a: x \mapsto \zeta_{2n}^a$ for some $a$, where $\zeta_{2n}$ is a primitive $2n$-th root of unity.
\[
V_{1, 1} \otimes V_{x^i, \rho_b} = V_{1, t_{-1}} \otimes V_{x^i, \rho_b}= V_{x^i, \rho_b} \quad \text{and} \quad
V_{1, t_{i}} \otimes V_{x^i, \rho_b} = V_{1, t_{-i}} \otimes V_{x^i, \rho_b}= V_{x^i, \rho_{m+b}}
\]
When $\chi$ is a two-dimensional irreducible character of $G$, we have
\[
V_{1, \chi_a} \otimes V_{x^i, \rho_b} = V_{x^i, \rho_{a+b}} \oplus V_{x^i, \rho_{a-b}}
\]
\end{proof}

\begin{Lemma}
The type-2-type-2 fusion rules are multiplicity free.
\end{Lemma}
\begin{proof}
The proof is analogous to that of Lemma \ref{dihedralevent2t2}
\end{proof}

\begin{Lemma}\label{dicevent1t3}
The type-1-type-3 fusion rules are multiplicity free.
\end{Lemma}
\begin{proof}
The proof is analogous to that of Lemma \ref{event1t3}.
\end{proof}

\begin{Lemma}
The type-2-type-3 fusion rules are multiplicity free.
\end{Lemma}
\begin{proof}
The proof is analogous to that of Lemma \ref{dihedralevent2t3}. We will only list the fusion rules here.

Let $\psi$ be a character of $C(g_L)$ which can only be one dimensional. Among the characters of $C(g_J)$, there are two choices of $\varphi$ such that $\psi(x^m)\varphi(x^m) = -1$. Call them $\varphi_1$ and $\varphi_2$. Then
\[
V_{x^i, \chi} \otimes V_{y, \psi} = V_{yx, \varphi_1} \otimes V_{yx, \varphi_2}.
\]
Now suppose $i$ is even. Then
\[
V_{x^i, \chi} \otimes V_{y, \psi} = V_{y, \varphi_1} \otimes V_{y, \varphi_2}.
\]
where $\varphi_1, \varphi_2$ are characters of $\{1, y, x^n, yx^n\}$ such that $\psi(x^n)\varphi_1(x^n) = \psi(x^n)\varphi_2(x^n) = -1$.
\end{proof}

\begin{Lemma}
The type-3-type-3 fusion rules are multiplicity free.
\end{Lemma}
\begin{proof}
The proof is analogous to that of Lemma \ref{oddt3t3}.
\end{proof}

Through these series of lemmas, we see that the fusion rules for $Q_{4n}$ when $n$ is odd are very similar to those of $D_{4n}$, but slightly simpler. In particular, we have proved

\begin{Theorem}
Let $G$ be a Dicyclic group $Q_{4n}$ where $n$ is odd. Then the category of $D(G)$-modules is multiplicity free.
\end{Theorem}

\section{Further Observations}

\subsection{Computations using Sage}
Our computations of various examples are facilitated by computer programming: we wrote Sage code to compute the $S$-matrix of $\mathrm{Rep}_f(D(G))$ by using \eqref{DGs} in Theorem \ref{DGMTC}, and implemented  two algorithms for computing the fusion coefficients: the first is by using the Verlinde formula, and the second is by using the character theoretic formula in Theorem \ref{generalfusion}. The first version of this code can be found \hyperlink{https://github.com/enqiLi/drinfeld_doubles}{here}\footnote{\url{https://github.com/enqiLi/drinfeld_doubles}}. The code was later edited to conform to the SageMath software programming conventions, and is expected to be included in SageMath version $10.1$ and later. The edited version can be found \hyperlink{https://github.com/sagemath/sage/pull/35387}{here}\footnote{\url{https://github.com/sagemath/sage/pull/35387}}.

\subsection{Other multiplicity free examples}
We attempted to find more groups whose Drinfeld doubles have multiplicity free representation categories. It is easy to see that Drinfeld doubles over abelian groups have multiplicity free representation categories. \cite{fusionrule} proved that extraspecial $2$-groups have the property, and we computed and found that all groups of order 16 have this property too. However, it is not true that all $2$-groups has this property: a counterexample with order 32 exists.

\subsection{Connection to the type B level 2 Verlinde algebras}

The type B Lie algebra is $\fk{g} = \fk{so}(2r+1)$ for some non-negative integer $r$. The Verlinde algebras of the level 2 representation category of type B Lie algebras arise as the fusion algebras of Wess-Zumino-Witten conformal field theories. These fusion rules are known, and can be found in \cite{NaiduRowell}. The following notations and results are directly from \cite{NaiduRowell}. Using the standard labeling of weights, let
\[
\lambda_1 = (1, 0, \cdots, 0), \cdots, \lambda_{r-1} = (1, 1, \cdots, 1, 0), \text{ and } \lambda_r = \frac{1}{2}(1, \cdots, 1).
\]
Then the representation category of $\fk{g}$ of level 2 has simple objects
\[
I = Y_0, Z = Y_{2\lambda_1}, Y_{\lambda_1}, \cdots, Y_{\lambda_{r-1}}, Y_{\lambda_{2r}}, X = Y_{\lambda_r}, X' = Y_{\lambda_r + \lambda_1}.
\]
The subscripts are the weights. For convenience, denote $Y_i = Y_{\lambda_i}$ for $1 \leq i \leq r-1$, and $Y_r = Y_{\lambda_{2r}}$. The fusion rules in this category resemble the fusion rules of the odd Dihedral groups to a remarkable extent: The objects $I$ and $Z$ behave like the one-dimensional type-1 objects in the Dihedral category, the $Y$'s behave like the type-2 objects and the two-dimensional type-1 objects, and $X$ and $X'$ behave like the type-3 objects. Here, the ``type 3'' fusion rules are (\cite{NaiduRowell}):
\begin{align*}
X \otimes X &= I \oplus \bigoplus_{i=1}^r Y_{i} \\
X \otimes X' &= Z \oplus \bigoplus_{i=1}^r Y_{i} \\
X \otimes Y_i &= X \oplus X' \\
X \otimes Z &= X'
\end{align*}
These fusion rules are exactly same as the ones in the odd Dihedral category. The fusion rules involving only $Y$'s also have the same form (i.e. type of objects and number of terms) as the ones in the odd Dihedral category, but they are not the same, and the fusion rings are not isomorphic. The striking similarity between these two fusion rings and the fact that they are in fact not isomorphic seem to suggest that they are part of a family of fusion rings with the fusion rules listed above, in particular with ``type-3'' objects giving fusion rules that involve almost all the simple objects.

\section*{Acknowledgements}

This paper is a condensed version of the author's undergraduate honors thesis at Stanford University under the supervision of Professor Daniel Bump. I would like to thank Professor Daniel Bump for providing numerous insights, suggestions, and guidance on writing.

\newpage

\printbibliography
\end{document}